\documentclass[preprints,article,accept,moreauthors,pdftex]{Definitions/mdpi}

\usepackage{amsmath}
\usepackage{times}
\usepackage{graphicx}
\usepackage{color}
\usepackage{multirow}
\usepackage{rotating}
\usepackage{bbm}
\usepackage{latexsym}
\usepackage{ulem}

\newcommand{\mbf}[1]{\mathbf{#1}}

\firstpage{1} 
\pubvolume{xx}
\issuenum{1}
\articlenumber{5}
\pubyear{2019}
\copyrightyear{2019}
\history{ }

\Title{The role of spectral complexity in connectivity estimation}

\Author{Elisabetta Vallarino $^{1}$, Alberto Sorrentino$^{1,2,}$, Michele Piana $^{1,2}$, and Sara Sommariva $^{1,}$}

\AuthorNames{Elisabetta Vallarino, Michele Piana, Alberto Sorrentino and Sara Sommariva}

\address{%
$^{1}$ \quad Dipartimento di Matematica, Università di Genova, Genova, Italy\\
$^{2}$ \quad CNR–SPIN, Genova, Italy}

\abstract{The study of functional connectivity from magnetoecenphalographic (MEG) data consists in quantifying the statistical dependencies among time series describing the activity of different neural sources from the magnetic field recorded outside the scalp. This problem can be addressed by utilizing connectivity measures whose computation in the frequency domain often relies on the evaluation of the cross-power spectrum of the neural time-series estimated by solving the MEG inverse problem. Recent studies have focused on the optimal determination of the cross-power spectrum in the framework of regularization theory for ill-posed inverse problems, providing indications that, rather surprisingly, the regularization process that leads to the optimal estimate of the neural activity does not lead to the optimal estimate of the corresponding functional connectivity. Along these lines, the present paper utilizes synthetic time series simulating the neural activity recorded by an MEG device to show that the regularization of the cross-power spectrum is significantly correlated with the signal-to-noise ratio of the measurements and that, as a consequence, this regularization correspondingly depends on the spectral complexity of the neural activity.}

\keyword{regularization theory, multivariate stochastic processes,
cross-power spectrum, magnetoencephalography, MEG,
functional connectivity, spectral complexity}

\begin{document}

\section{Introduction}
Magnetoencephalography (MEG) provides high temporal resolution measurements of the magnetic field associated to neural currents. The MEG device relies on superconducting sensors, named SQUIDs, organized in a helmet array close and around the scalp. MEG experimental time series can be used essentially to address two neuroscientific problems, whose solution requires both an accurate mathematical modelization based on Maxwell's equations, and the numerical reduction of such formal models \citep{hamalainen1993magnetoencephalography}. 

The first problem is concerned with the dynamical ill-posed inverse problem of estimating parameters associated to the neural sources inducing the magnetic field signal \citep{hail94,calvetti2015hierarchical, costa2017bayesian, sorrentino2017inverse,bekhti2018hierarchical, luria19, ilsa19}. The second problem is concerned with the quantification of the interactions among neural sources located in different cortical areas and intertwined by means of either anatomical or functional connectivity \citep{geweke1982measurement,baccala2001partial, nolte04, pereda_05, sakkalis_11,chella2014third}.

In particular, the connectivity problem can be addressed by either computing proper connectivity metrics directly from the experimental time series provided by the MEG sensors or searching for connections in the source space, i.e. among the neural time series estimated as solutions of the inversion process. This second approach has the advantages of reducing the impact of volume conduction and providing results that can be more easily interpreted in the framework of neuroscientific models \citep{schoffelen2009source,barzegaran2017functional, van2019critical}. Several approaches for identifying connectivity paths rely on physiological models assuming that the functional communication between different brain areas is regulated by the synchronization of their activity at specific temporal frequencies \cite{fries2005, fries2015}. This implies that, for these models, the frequency domain represents the natural computational framework where to perform the connectivity analysis. This is the reason why, in the present paper, we focus on the analysis of the cross-power spectrum, which is the mathematical quantity of reference for the computation of most frequency-domain connectivity measures \citep{nunez99,nolte04,nolte08}. From an operational viewpoint, the computation of the cross-power spectrum in the source space typically relies on a two-step procedure: first the neural activity is estimated by applying a regularized inversion method on the recorded time series and then the cross-power spectrum is computed from the Fourier transform of the estimated neural time series \citep{sc_gr19}.  

\begin{figure} 
    \centering
    \includegraphics[scale=1]{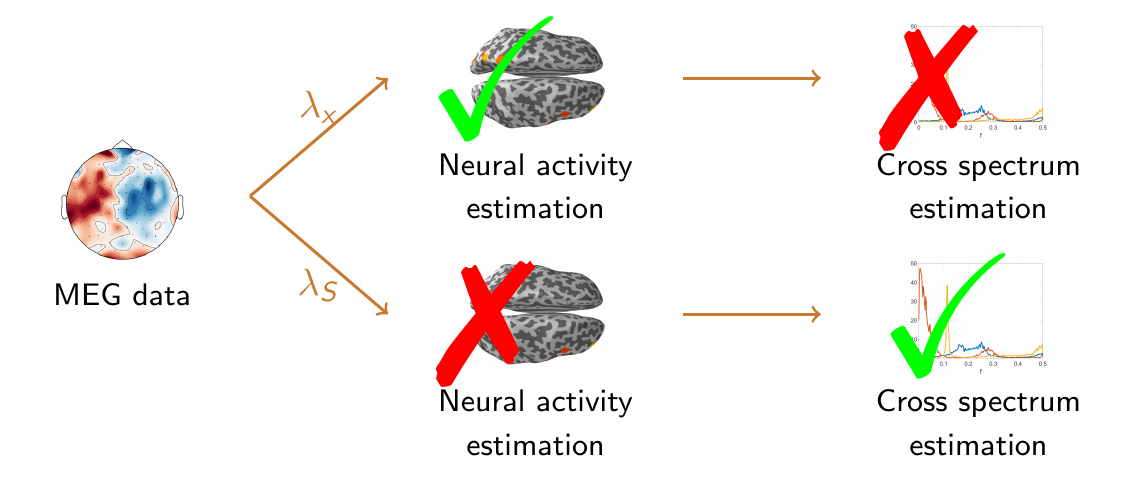}
    \caption{{Schematic representation of the differences between the regularization parameter providing the best time series estimate ($\lambda_\mbf{x}$) and the one providing the best cross-power spectrum estimate ($\lambda_\mbf{S}$). The first one provides an optimal reconstruction of the neural activity, but it may not lead to an optimal estimate of the cross-power spectrum; vice versa $\lambda_\mbf{S}$ provides an optimal reconstruction of the cross-power spectrum at the expense of a sub-optimal estimate of the time series.} }\label{fig:two_step}
\end{figure}

This paper investigates how to optimize the inversion procedure in order to obtain the best possible estimate of the neural cross-power spectrum. In fact, we consider the Tikhonov method (better known as Minimum Norm Estimation (MNE) in the MEG world \citep{hail94}) as a paradigmatic inversion technique and we study the interplay between the regularization parameter providing the reconstructed neural time series minimizing the relative error in $\ell_2$-norm, and the one that allows the optimal estimate of the cross-power spectrum according to the normalized Frobenius norm. The conceptual motivation of this problem is illustrated in Figure \ref{fig:two_step}, which tentatively sketches the result of recent investigations in MEG-based connectivity research, i.e. that the regularization parameter leading to the optimal estimate of the neural activity may not lead to the optimal estimate of the cross-power spectrum and, vice versa. In fact, in \cite{hietal16} the authors used numerical simulations to compare the parameter that provides the best estimate of the power spectrum with the one that provides the best estimate of coherence and showed that the latter is in general two orders of magnitude smaller than the former. More recently, \citet{vallarino2020two} addressed an analogous problem via analytical computations, considering a simplified model. Specifically, under the assumption that the neural time series are realizations of white Gaussian processes, the authors proved that the parameter providing the best neural activity estimate is more than twice as large as the one providing the best estimate of the cross-power spectrum. 

The present paper focuses on an analysis of the impact of spectral complexity of the actual neural signal on the value of the two regularization parameters. Specifically, we simulate synthetic MEG signals and discuss how the optimal parameter for the reconstruction of the cross-power spectrum depends on its signal-to-noise ratio and how this latter quantity is related to the spectral richness of the neural sources. To this aim, we considered a simulation setting in which the signal is modelled as a multivariate autoregressive process. 

The plan of the paper is as follows. Section 2 introduces the problem in a formal way. Section 3 describes how the synthetic data are simulated and analysed. Section 4 presents the results of the analysis. Our conclusions are offered in Section 5.

\section{Definition of the problem}

\subsection{Forward model}

Let $\mbf{X}(t)=(X_1(t),\dots,X_N(t))^\top\in{R}^{N}$ be a multivariate stationary stochastic process whose realisations $\mbf{x}(t)$ can not be observed and let $\mbf{Y}(t)=(Y_1(t),\dots,Y_M(t))^\top\in {R}^{M}$ be the process whose realizations $\mbf{y}(t)$ are used to infer information on $\mbf{x}(t)$. Let $\mbf{Y}(t)$ and $\mbf{X}(t)$ be related by the following equation
\begin{equation} \label{eq:fwd_mod}
    \mbf{Y}(t) = \mbf{G}\mbf{X}(t)+\mbf{N}(t) ~~,
\end{equation}
where $\mbf{G}\in {R}^{M\times N}$ is the forward matrix and $\mbf{N}(t)=(N_1(t),\dots,N_M(t))^\top\in {R}^{M}$ is the measurement noise, that is here assumed to be a white Gaussian process with zero mean and covariance matrix $\alpha^2\mbf{I}$, i.e. $\mbf{N}(t)\sim\mathcal{N}(0,\alpha^2\mbf{I})$, independent from $\mbf{X}(t)$.

\subsection{Cross-power spectrum}

We are interested in reconstructing the cross-power spectrum of $\mbf{X}(t)$, which describes the statistical dependencies between each pair of time series  $(X_j(t),X_k(t))_{j,k \in\{1,\dots,N\}}$. The cross-power spectrum is a one parameter family of $N\times N$ matrices $\mbf{S}^\mbf{X}(f)$, whose $(j,k)-th$ element is defined as
\begin{equation}
  S^\mbf{X}_{j,k}(f) = \lim_{T\to\infty} \frac{1}{T} E[\hat{X}_j(f,T)\hat{X}_k(f,T)^H],
\end{equation}
where $\hat{X}_j(f,T)$ is the Fourier transform of $X_j(t)$ over the interval $[0,T]$, defined as
\begin{equation}
    \hat{X}_j(f,T) = \int_{0}^{T} X_j(t) e^{-2\pi ift} \textup{d}t
\end{equation}
and $X^H$ is the Hermitian transpose of $X$ \citep{be_pi11}

{Given a realization $\mbf{x}(t)$ of the process $\mbf{X}(t)$, the cross-power spectrum $S^\mbf{X}(f)$ can be estimated via the Welch method \citep{welch1967use}, which consists in partitioning the data in $P$ overlapping segments multiplied by a window function, $\{w(t)\mbf{x}^p(t)\}_{p=1}^P$, computing their discrete Fourier transform $\hat{\mbf{x}}^p(f) = \frac{1}{L}\sum_{t=0}^{L-1}\mbf{x}^p(t)w(t)e^{\frac{ -2\pi itf}{L}}$ and averaging: 
    \begin{equation}\label{S_lam}
        \mbf{S}^{\mbf{x}}(f) = \frac{L}{PW}\sum_{p=1}^P \mbf{\hat{x}}^p(f)\mbf{\hat{x}}^p(f)^H, ~~~~ f=0,\dots,L-1,
        \end{equation}
    where L is the length of each segment and $W=\frac{1}{L} \sum_{t=0}^{L-1}w(t)^2$.}

It is often the case that the data reaches high dimension, and visual inspection of the cross-power spectrum is not doable. In such cases a metric that describes the spectral properties of the signals would be useful. Here we use the spectral complexity coefficient, defined as follows.\\

\begin{Definition}\label{def:sc_coeff}
Given a realization $\mbf{x}(t)$ of the process $\mbf{X}(t)$, and the corresponding cross-power spectrum $\mbf{S}^{\mbf{x}}(f)$, we define the spectral complexity coefficient as the average of the elements of the upper triangular part of the matrix obtained by computing the squared $\ell_2-$norm over the frequencies of  $\mbf{S}_{j,k}^\mbf{x}(f)$, $j,k=1, \dots, N$, that is
\begin{equation}
    c = \frac{2}{N(N+1)}\sum_{j=1}^N\sum_{k=j}^N\sum_{f} \left| S_{j,k}^\mbf{x}(f) \right|^2.
\end{equation}
\end{Definition}

 The spectral complexity coefficient assumes small values if the elements of the cross-power spectrum are flat, that is when time series do not present any periodic trend and no dependencies among the pairs of time series are present. On the contrary, it assumes large values if the elements of the cross-power spectrum are peaked, that is when time series present periodic trends and complex relations among them.
Finally, we observe that in Definition \ref{def:sc_coeff} only the elements on the upper triangular part of $\mbf{S}^\mbf{x}(f)$ are considered because $\mbf{S}^\mbf{x}(f)$ is Hermitian.

\subsection{Two-step approach for cross-power spectrum estimation}
Let us now consider a realization of equation (\ref{eq:fwd_mod}). Further than an estimate of the hidden data $\mbf{x}(t)$, an estimate of the cross-power spectrum can be obtained from $\mbf{y}(t)$. Such estimate can be achieved through a two-step process \citep{sc_gr19}:
\begin{itemize}
    \item[i.] First, a regularized estimate $\mbf{
    x}_\lambda(t)$ of $\mbf{x}(t)$ is obtained by solving the inverse problem associated to equation (\ref{eq:fwd_mod}). Here we consider the Tikhonov regularized solution \citep{ti_etal13} of the problem which is defined as
    \begin{equation} \label{eq:x_lam}
        \mbf{x}_\lambda(t) = \text{arg}\,\min\limits_{\mbf{x}(t)}\left\{ \left\Vert \mbf{Gx}(t)-\mbf{y}(t)\right\Vert_2^2+\lambda\left\Vert \mbf{x}(t) \right\Vert_2^2 \right\};
    \end{equation}
    where $\lambda$ is a proper regularization parameter and $\Vert \cdot\Vert_2$ is the $\ell_2$-norm.
    \item[ii.] Then, the corresponding estimate of the cross-power spectrum $\mbf{S}^{\mbf{x}_\lambda}(f)$ is computed from the reconstructed time series using the Welch method, as described in the previous section.
    
\end{itemize}

When applying this two-step process the regularization parameter $\lambda$ in equation (\ref{eq:x_lam}) has to be set for the computation of $\mbf{x}_\lambda(t)$. Thus, the problem naturally arises of the choice of such parameter, which can be set in order to optimally reconstruct either $\mbf{x}_\lambda(t)$ or $\mbf{S}^{\mbf{x}_\lambda}(f)$. We define optimality through the minimization of the normalised norm of the discrepancy between the true and the reconstructed time series and cross-power spectra as follows.

\begin{Definition} \label{def:optimality}
Given the regularized solution (\ref{eq:x_lam}) and the cross-power spectrum (\ref{S_lam}), we define the optimal regularization parameter for the reconstruction of $\mbf{x}(t)$ as
\begin{equation} \label{eq:err_x}
    \lambda_{\mbf{x}}^* = \arg \min_{\lambda}\varepsilon_{\mbf{x}}(\lambda) ~~~ \textrm{with} ~~~  \varepsilon_{\mbf{x}}(\lambda) = \frac{\sum_t \left\| \mbf{x}_{\lambda}(t) - \mbf{x}(t) \right\|_2^2}{\sum_t \left\| \mbf{x}_{\lambda}(t)\right\|_2^2+\sum_t\left\|\mbf{x}(t) \right\|_2^2} ~~~;
\end{equation}
and the optimal parameter for the reconstruction of $\mbf{S}^{\mbf{x}}(f)$ as
\begin{equation}\label{eq:err_S_x}
    \lambda_{\mbf{S}}^* = \arg \min_{\lambda} \varepsilon_{\mbf{S}}(\lambda) ~~~ \textrm{with} ~~~ \varepsilon_{\mathbf{S}}(\lambda) = \frac{\sum_f \left\| \mbf{S}^{\mbf{x}_{\lambda}}(f) - \mbf{S}^{\mbf{x}}(f) \right\|_F^2}{\sum_f \left\| \mbf{S}^{\mbf{x}_{\lambda}}(f)\right\|_F^2+\sum_f \left\| \mbf{S}^{\mbf{x}}(f) \right\|_F^2}~~~;
\end{equation}
where {$\|\cdot\|_F$ is the Frobenius norm;} $\varepsilon_{\mbf{x}}(\lambda)$ and $\varepsilon_{\mbf{S}}(\lambda)$ will be called reconstruction errors.

\end{Definition}

The reconstruction errors range from 0 to 1 and penalize both a too small and a too large value of $\lambda$. In fact, they assume their maximum value when either $\lambda$ is very high and thus $\mbf{x}_\lambda(t)$ is negligible with respect to $\mbf{x}(t)$, or when $\lambda$ is too small and thus, vice versa, $\mbf{x}(t)$ is negligible with respect to $\mbf{x}_\lambda(t)$.
This definition may appear overly complex compared to, e.g., a mere $\ell_2$-norm of the difference; however, in the presence of sparse data where only few time series are non-zero, the simple $\ell_2$-norm would prefer a very high regularization parameter in order to minimize the error on the null time series, at the expense of the error on the non-zero ones; our definition aims to cope with this limitation of the $\ell_2$-norm. A similar definition has been introduced in \cite{chella2019impact}.

In experimental contexts, where $\mbf{x}(t)$ is not known, the choice of the optimal regularization parameter is crucial. This matter is widely discussed in literature \citep{thompson1991study, hanke1993regularization, hansen1998rank, vogel2002computational}, and many criteria have been proposed. Such criteria apply to equation (\ref{eq:fwd_mod}) and can be used to set the regularization parameter $\lambda_\mbf{x}$. A possibility is to set the regularization parameter as a function of the signal-to-noise ratio (SNR), which describes the level of the desired signal with respect to that of the measurement noise; for equation (\ref{eq:fwd_mod}) the SNR is defined as follows. 

\begin{Definition}
Consider the linear model (\ref{eq:fwd_mod}). We define the signal-to-noise ratio of $\mbf{X}(t)$ related to such model as
\begin{equation} \label{eq:SNR_X}
    \textup{SNR}^{\mbf{X}} = 10\log_{10}\left(\frac{\sum_{t}\left\|\mbf{G}\mbf{X}(t)\right\|_2^2}{\sum_{t}\left\|\mbf{N}(t)\right\|_2^2}\right).
\end{equation}
\end{Definition}

To the best of our knowledge, the choice of the optimal regularization parameter for the reconstruction of the cross-power spectrum has never been related to the signal-to-noise ratio. This relation will be presented in Section \ref{section:results}; however we first need to relate the cross-power spectrum of the unknown $\mbf{S}^\mbf{X}(f)$ with that of the data $\mbf{S}^\mbf{Y}(f)$. 

By computing the cross-power spectrum of both sides of equation (\ref{eq:fwd_mod}) and from the linearity of the Fourier transform it follows that
\begin{equation} \label{eq:specrum_model}
    \mbf{S}^{\mbf{Y}}(f) = \mbf{G}\mbf{S}^{\mbf{X}}(f)\mbf{G}^\top + \mbf{S}^{\mbf{N}}(f),
\end{equation}
where the mixed terms $\mbf{S}^\mbf{XN}(f)$ and $\mbf{S}^\mbf{NX}(f)$ are negligible thanks to the independence between $\mbf{X}(t)$ and $\mbf{N}(t)$. Just like for equation (\ref{eq:fwd_mod}), we can define the signal-to-noise ratio for equation (\ref{eq:specrum_model}) as follows.

\begin{Definition}
Consider the linear model (\ref{eq:specrum_model}). We define the signal-to-noise ratio of $\mbf{S}^{\mbf{X}}(f)$ related to such model as
\begin{equation}\label{eq:SNR_S}
    \textup{SNR}^{\mbf{S}} = 10\log_{10}\left(\frac{\sum_{f}\left\|\mbf{G}\mbf{S}^{\mbf{X}}(f)\mbf{G}^\top\right\|_F^2}{\sum_{f}\left\|\mbf{S}^{\mbf{N}}(f)\right\|_F^2}\right).
\end{equation}
\end{Definition}

 This definition is in line with the definition of $\textup{SNR}^\mbf{X}$ for the signal, the main difference being in the use of the Frobenius norm rather than the $\ell_2$-norm, motivated by the fact that we are working with matrices rather than vectors.
\section{Generation and analysis pipeline of the MEG simulated data.} \label{sec:simulations}

In this section we will describe the numerical simulation that led to the main results of our study. First we introduce the continuous MEG forward problem and its discretized version, then we describe how we generated the data and, finally, we describe the inverse model and how we numerically computed the optimal regularization parameters.

\subsection{MEG forward model.}

The MEG forward problem aims at computing the magnetic field produced outside the head by an electric current that flows inside the brain. The quasi static approximation of Maxwell's equations provides the local relationship between the recorded magnetic field and the neural currents \cite{hamalainen1993magnetoencephalography, baillet01, sorrentino10}. The two equations that are of interest here read as

\begin{eqnarray}
& \nabla \times\mbf{E}(\mbf{r},t) =0\label{eq:maxwell_2}\\ 
&\nabla\times\mbf{B}(\mbf{r},t) = \mu_0\mbf{J}(\mbf{r},t);\label{eq:maxwell_4}
\end{eqnarray}

where $\mbf{E}(\mbf{r},t)$ and $\mbf{B}(\mbf{r},t)$ are the electric and magnetic fields at location $\mbf{r}$ and time $t$, $\mu_0$ is the magnetic permeability in vacuum and $\mbf{J}(\mbf{r},t)$ is the total electric current that flows inside the brain. The latter is the sum of two contributions 
\begin{equation}
\mbf{J}(\mbf{r},t) = \mbf{J}^p(\mbf{r},t)+ \mbf{J}^v(\mbf{r},t), 
\end{equation}
$\mbf{J}^p(\mbf{r},t)$ being the primary current directly related to the brain activity, while $\mbf{J}^v(\mbf{r},t) = -\sigma(\mbf{r})\nabla V(\mbf{r},t)$ is the induced volume current due the non-null conductivity $\sigma(\mbf{r})$ of the brain, $V(\mbf{r},t)$ being the electric scalar potential.

The manipulation of Maxwell's Equations leads to the Biot-Savart Equation
\begin{equation}
\mbf{B}(\mbf{r},t) =\mbf{B}_0(\mbf{r},t)- \frac{\mu_0}{4\pi}\int_{\Omega}\sigma(\mbf{r}')\nabla'V(\mbf{r}',t)\times \frac{\mbf{r}-\mbf{r}'}{|r-r'|^3}\textup{d}v',
\end{equation}
where $\Omega$ is the volume occupied by the brain, the first term $\mbf{B}_0(\mbf{r},t)=\frac{\mu_0}{4\pi}\int_{\Omega}\mbf{J}(\mbf{r}',t)\frac{\mbf{r}-\mbf{r}'}{|r-r'|^3}\textup{d}v'$ is the magnetic field {induced by} the primary current while the second term is {related to} the volume current.

Solving the forward problem requires the computation of these two contributions knowing the primary current. While for the first one straightforward numerical integration is feasible, for the second one it is common to model the head as the union of nested homogeneous volumes $\{\Omega_j\}_{j=1,\dots,J}$ and to replace volume integration with surface integration. In this way Biot-Savart equation becomes
\begin{equation}\label{eq:B_final}
    \mbf{B}(\mbf{r},t) = \mbf{B}_0(\mbf{r},t)+\frac{\mu_0}{4\pi} \sum_{i,j}(\sigma_i-\sigma_j)\int_{\partial \Omega_{i,j}} V(\mbf{r}',t) \frac{\mbf{r}-\mbf{r}'}{|r-r'|^3}\times \mbf{n}_{i,j}(\mbf{r}')\textup{d}s',
\end{equation}
where $\partial \Omega_{i,j}$ is the contact surface between regions $\Omega_i$ and $\Omega_j$ and $\mbf{n}_{i,j}(\mbf{r}')$ is the unit vector normal to the surface $\partial \Omega_{i,j}$ at $\mbf{r}'$ from
region $i$ to region $j$.

The forward problem can now be solved by computing the second term at the right hand side of equation (\ref{eq:B_final}) after having computed $V(\mbf{r},t)$ by solving the equation
\begin{equation}\label{locale}
\nabla\cdot\mbf{J}^p(\mbf{r},t)-\nabla\cdot\left(\sigma(\mbf{r})\nabla V(\mbf{r},t)\right)=0~,
\end{equation}
which follows from equation (\ref{eq:maxwell_4}) by applying the divergence.

For further details on the MEG forward problem we refer the reader to \cite{hamalainen1993magnetoencephalography}.

\subsection{The leadfield matrix}
Experimental contexts require the discretization of the forward problem. This involves a discretization of both the volume occupied by the brain and the volume outside the head.

When using a distributed model for the primary current $\mathbf{J}^p$, the brain volume is uniformly divided in $N$ small parcels. If $N$ is sufficiently big and thus each parcel has a sufficient small area, the activity in each brain parcel is approximated by a point-like source, henceforth denoted as dipole. From a mathematical point of view, each dipole is a vector whose strength and direction represent the intensity and orientation of the primary current in the corresponding brain area \cite{hamalainen1993magnetoencephalography}.

As for the volume outside the brain, it is natural to discretize it in correspondence of the MEG sensors. Let us denote the measured magnetic field as $\mbf{y}(t)=(y_1(t),\dots,y_M(t))$. Now, observing that the magnetic field $\mbf{B}$ depends linearly on the primary current $\mbf{J}^p$, the magnetic field in correspondence of the sensors of the instrument is 
\begin{equation}\label{eq:MEG_model1}
    \mbf{y}(t) = \sum_{k=1}^N{G}(\mbf{r}_k)\mbf{q}_k(t) + \mbf{n}(t),
\end{equation}
where $\mbf{r}_k$, $k=1, \dots, N$, is the location of the $k$-th brain parcel, ${G}(\mbf{r}_k)\in R^{M\times3}$ is the corresponding leadfield matrix, and $\{\mbf{q}_k(t)\}_{k=1,\dots,N}$ are the electric current intensities along the three orthogonal direction of the $N$ dipoles within the brain at time $t$  and $\mbf{n}(t)$ is the measurement noise. The $l$-th column of ${G}(\mbf{r}_k)$ contains the measurement at a sensor level when a unit current dipole is placed at location $\mbf{r}_{k}$ and oriented along the $l$-th orthogonal direction.

In this work we assume dipoles to be located only on the brain cortical mantle and their orientation to be normal to the local cortical surface \cite{dale93}. In this case the electric current intensities are scalar quantities (we refer to them as $\{q_k\}_{k=1,\dots,N}$) and the leadfield matrices are column vectors (we refer to them as $\{G_k\}_{k=1,\dots,N}$).

Let us define
\begin{equation} \label{eq:x}
    \mbf{x}(t) := (q_1(t),\dots, q_N(t))
\end{equation}
and
\begin{equation} \label{eq:G}
    \mbf{G} := [G_1,\dots, G_N] \in R^{M\times N}~;
\end{equation}
reassembling equations (\ref{eq:x}) and (\ref{eq:G}) in to equation (\ref{eq:MEG_model1}), we get 
\begin{equation}
    \mbf{y}(t) = \mbf{G}\mbf{x}(t) +\mbf{n}(t),
\end{equation}
which can be interpreted as a realization of equation (\ref{eq:fwd_mod}). From now on we will refer to $\mbf{G}$ as to the leadfield matrix.\\

For the simulation presented in this work we used the leadfield matrix available in the sample dataset of MNE Python \citep{gr_etal14}. We selected magnetometers and set a fixed orientation. For computational reasons, the available source space, containing $1884$ sources, was uniformly down-sampled to obtain $274$ sources. Thus, our model has $M=102$ sensors and $N=274$ dipole sources.

\subsection{Data generation}
\begin{figure} 
    \centering
    \includegraphics[scale=0.3]{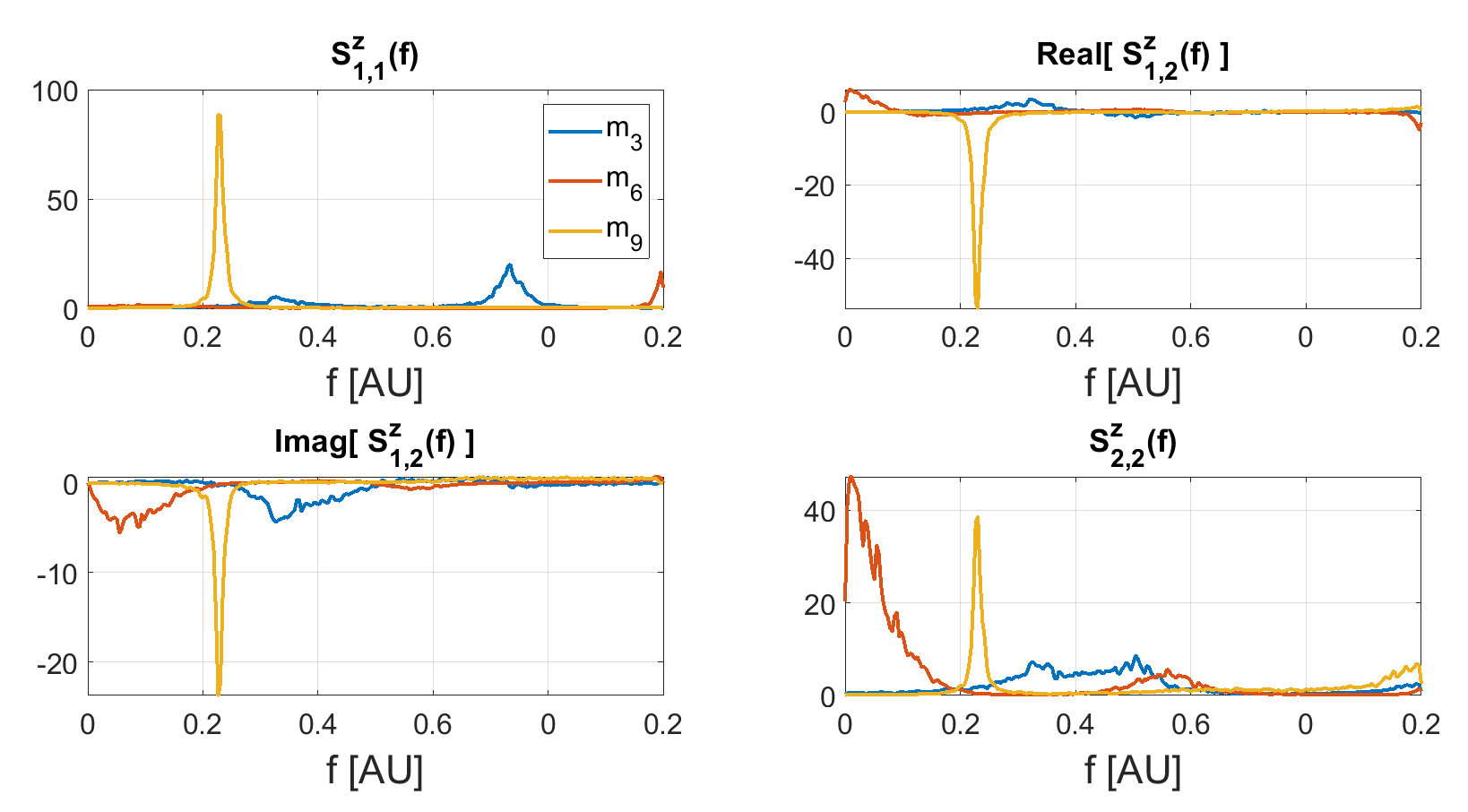}
    \caption{Real and imaginary part of the cross-power spectra of three simulated time series. Higher values of spectral complexity correspond to more peaked spectra.}
    \label{fig:cpsd}
\end{figure}

We simulated $N_{mod}=10$ pairs of active sources, $(z_1(t),z_2(t))^\top$, with unidirectional coupling from the first to the second; their time series follow a multivariate autoregressive (MVAR) model of order $P=5$ \citep{lutkepohl2005new, soetal19}
\begin{equation}
\left(\begin{array}{c}
    z_1(t) \\
    z_2(t) 
    \end{array} \right)
    =
    \sum_{k=1}^{P} \left( 
    \begin{array}{cc}
    a_{1,1}(k) & 0\\
    a_{2,1}(k) & a_{2,2}(k)
    \end{array} \right)
    \left(\begin{array}{c}
    z_1(t-k) \\
    z_2(t-k) 
    \end{array} \right)
    +
    \left(\begin{array}{c}
    \varepsilon_1(t) \\
    \varepsilon_2(t) 
    \end{array} \right), ~~ t=1,\dots,T.
\end{equation}
The non-zero elements $a_{i,j}(k)$ of the coefficient matrices were drawn from a normal distribution of  zero mean and standard deviation $\gamma$, and T=10000. We retained only coefficient matrices providing (i) a stable MVAR model \citep{lutkepohl2005new} and (ii) pairs of signals $(z_1(t),z_2(t))^\top$ such that the $\ell_2$-norm of the strongest one was less than 3 times the $\ell_2$-norm of the weakest one. In order to obtain time series with different spectral complexity coefficients we set $\gamma$ to $N_{mod}$ different values randomly drawn in the interval $[0.1, 1]$. The values of the spectral complexity coefficient of the $N_{mod}$ simulated time-series are reported in Table \ref{tab:coefficients}. Finally, the resulting time series $(z_1(t),z_2(t))^\top$ were normalized by the mean of their standard deviations {over} time, so that pairs of time series drawn from different models had similar magnitude. Figure \ref{fig:cpsd} shows a sample of the the cross-power spectra {among the simulated pairs of time series. The Figure shows that for increasing values of the spectral complexity coefficient the cross-power spectrum of the corresponding time series becomes more peaked.} Each pair of simulated time series was then assigned to $N_{loc} = 20$ pairs of point like sources randomly chosen in the source space, so that the ratio of the norms of the corresponding columns of the leadfield matrix was close to one, i.e. they had similar intensity at a sensor level, and their distance was grater than $7\ \textup{cm}$. The remaining $N-2$ sources were set to have null activity. 

Source space activity was then projected to sensor level by multiplying the simulated source activity by the leadfield matrix and white Gaussian noise was added to obtain $N_{snr}=6$ levels of $\textup{SNR}^\mbf{X}$ evenly spaced in the interval $[-20 \textup{dB}, 5\textup{dB}]$. 

Summarizing, we generated $N_{mod}\cdot N_{loc}\cdot N_{snr} = 1200$ different sensor level configurations. The green box in Figure \ref{fig:pipeline} shows a visual representation of the simulation pipeline.

\begin{table}[H]
\caption{The table reports the values of the spectral complexity coefficients, $c_j$, associated to each simulated MVAR model, $m_j$, $j = 1, \dots, N_{mod}$.}\label{tab:coefficients}
\centering
\begin{tabular}{c|*{10}c}
\toprule

\textbf{Model}		& 1 &2 &3 &4 &5 &6 &7 &8 &9 &10\\
\midrule
\textbf{Spectral complexity coefficient}		&1.41 &1.96 &2.14 &3.10 &3.44 &4.17 &4.64 &5.67 &6.69 &8.67\\
\bottomrule
\end{tabular}
\end{table}

\subsection{Inverse model}
\begin{figure} 
    \centering
    \includegraphics[scale=1.1]{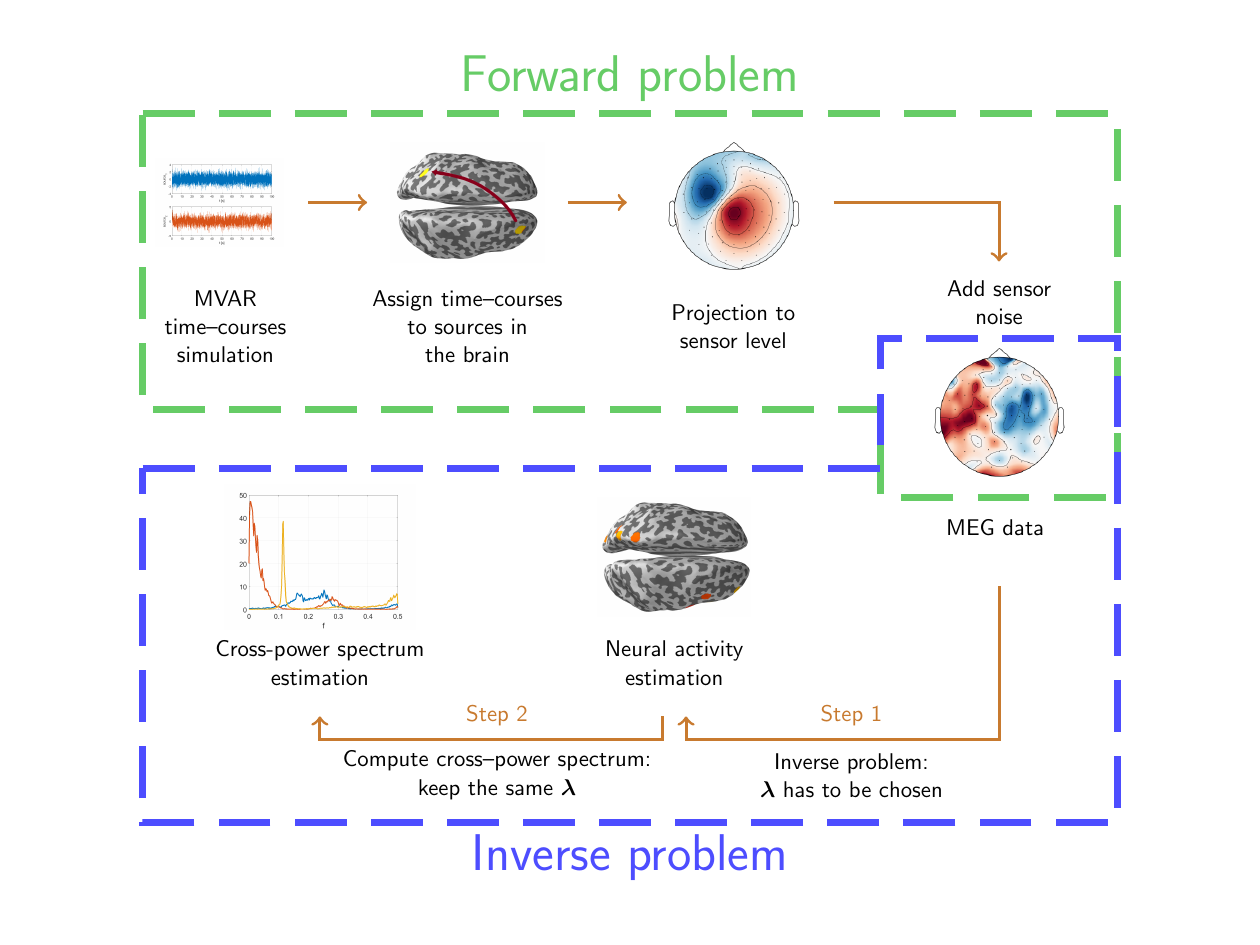}
    \caption{Pipeline of the simulation of the data (green box) and of the estimation of the cross-power spectrum (blue box).}\label{fig:pipeline}
\end{figure}

Source space time series were reconstructed using {the Tikhonov method, also known as} minimum norm estimate (MNE)  \citep{hail94} within the MEG community. For each combination of source time series, source locations and $\textup{SNR}^\mbf{X}$ level we computed the optimal regularization parameters $\lambda_\mbf{x}^*$ and $\lambda_\mbf{S}^*$ by minimizing the reconstruction errors $\varepsilon_\mbf{x}(\lambda)$ and $\varepsilon_\mbf{S}(\lambda)$, defined in Definition \ref{def:optimality}. The minimization procedure was achieved by using the Matlab built in function \texttt{fminsearch}. The blue box in figure \ref{fig:pipeline} describes the inverse procedure to obtain an estimate of the cross-power spectrum and stresses the role of the regularization parameter in the two-step process. 

\section{Results}\label{section:results}
In this section we illustrate the results of our analysis. We will begin with the description of the analytical dependence between $\textup{SNR}^\mbf{X}$ and $\textup{SNR}^\mbf{S}$, then we will highlight how the optimal parameter for the reconstruction of the cross-power spectrum depends on $\textup{SNR}^\mbf{S}$ and how this implies that the spectral complexity of the signal is behind such dependence. As a byproduct, this analysis also confirms the results of  \citet{vallarino2020two} in the case of a more complex setting.

\subsection{Analytical relation between $\textup{SNR}^\mbf{X}$} and $\textup{SNR}^\mbf{S}$
From equations (\ref{eq:SNR_X}) and (\ref{eq:SNR_S}) and reminding that $\mbf{N}(t)\sim \mathcal{N}(0,\alpha^2\mbf{I})$ it follows that
\begin{equation} \label{eq:SNR_X_alpha}
    \textup{SNR}^{\mbf{X}} = 10\log_{10}\left(\frac{\sum_{t}\left\|\mbf{G}\mbf{X}(t)\right\|_2^2}{MT\alpha^2}\right);
\end{equation}
and
\begin{equation}\label{eq:SNR_S_alpha}
    \textup{SNR}^{\mbf{S}} = 10\log_{10}\left(\frac{\sum_{f}\left\|\mbf{G}\mbf{S}^{\mbf{X}}(f)\mbf{G}^\top\right\|_F^2}{MN_f\alpha^2}\right),
\end{equation}
where $T$ is the number of time points and $N_f$ is the number of frequencies used to compute the cross-power spectrum. Observe that to derive equation (\ref{eq:SNR_S_alpha}) we used the fact that the cross-power spectrum of a white noise Gaussian process of zero mean and covariance matrix $\alpha^2 \mbf{I}$ is $\mbf{S}^\mbf{N}(f) = \alpha^2\mbf{I}$.

By isolating $\alpha^2$ from equation (\ref{eq:SNR_X_alpha}) and substituting in equation (\ref{eq:SNR_S_alpha}) we obtain
\begin{equation}\label{eq:SNR_relation}
    \textup{SNR}^{\mbf{S}} = 10\log_{10}\left( \frac{T^2 M \sum_{f}\left\| \mbf{G}\mbf{S}^\mbf{X}(f)\mbf{G}^\top\right\|_F^2}{N_f\sum_{t}\left\|\mbf{G}\mbf{X}(t)\right\|_2^4}\right)+2\textup{SNR}^\mbf{X}.
\end{equation}

Equation (\ref{eq:SNR_relation}) relates the the signal-to-noise ratio of $\mbf{X}(t)$ with that of $\mbf{S}^\mbf{X}(f)$. It shows that, for same levels of $\textup{SNR}^\mbf{X}$, $\textup{SNR}^\mbf{S}$ changes with the spectral complexity coefficient of the signals. In fact, the higher the spectral complexity coefficient, the higher the quantity $\left\| \mbf{G}\mbf{S}^\mbf{X}(f)\mbf{G}^\top\right\|^2_F$. Intuitively, this happens because when the signal has a higher spectral complexity coefficient its cross-power spectrum is more peaked and thus it is stronger over the cross-power spectrum of the noise with respect to a signal with a lower spectral complexity coefficient.

\subsection{Dependence of $\lambda_\mbf{S}^*$ on $\textup{SNR}^\mbf{S}$}

\begin{figure} 
    \centering
    \includegraphics[scale=0.3]{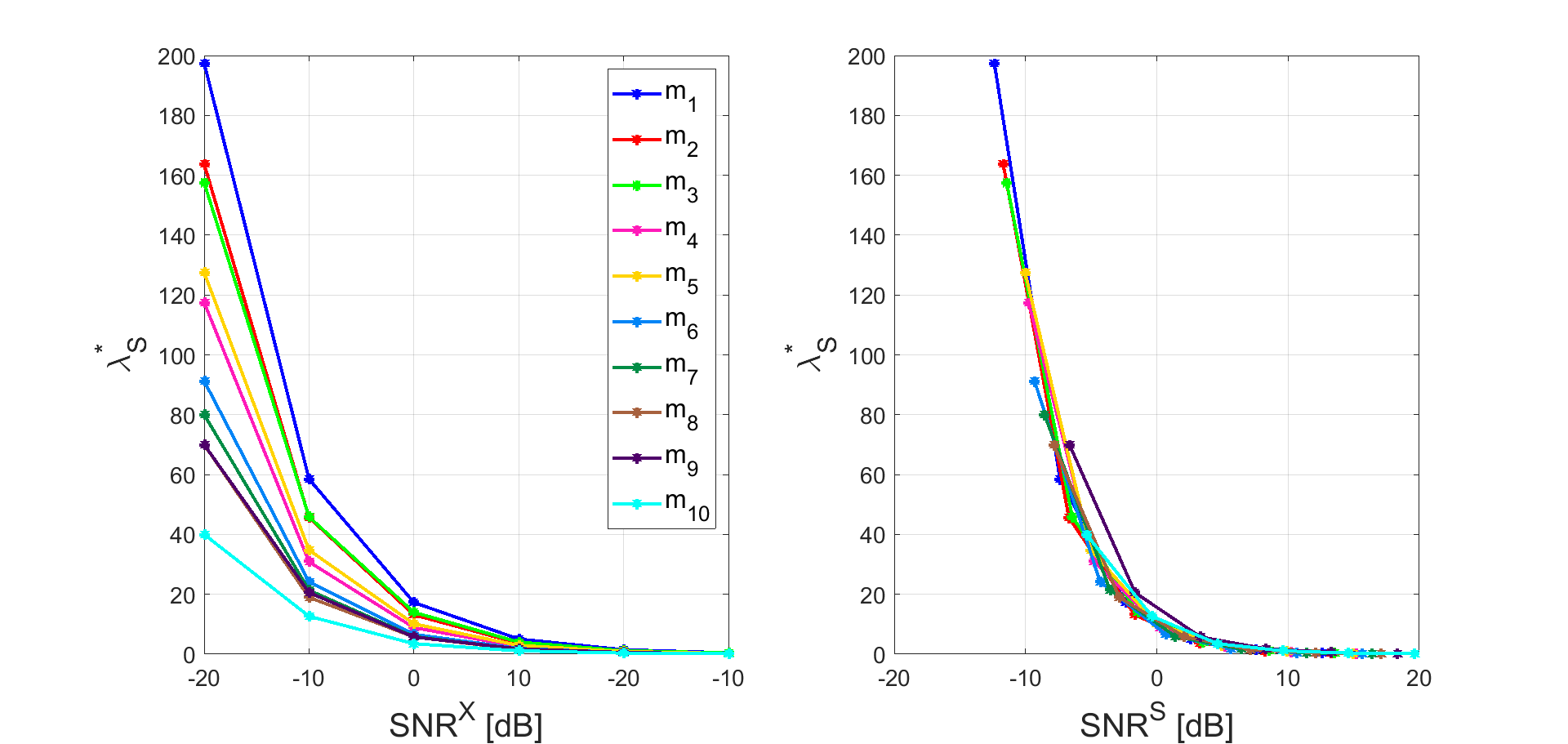}
\caption{Optimal regularization parameters for the reconstruction of the cross-power spectrum ($\lambda_\mbf{S}^*$) as a function of $\textup{SNR}^\mbf{X}$ (left) and $\textup{SNR}^\mbf{S}$ (right). Different colors correspond to different MVAR models. On the left, the lines have different heights, while on the right they overlap, meaning that the dependence of $\lambda_\mbf{S}^*$ on  $\textup{SNR}^\mbf{S}$ is neater with respect to $\textup{SNR}^\mbf{X}$.}
        \label{fig:lamS_SNR_S}
\end{figure}

As described in Section \ref{sec:simulations} we simulated several sensor level configurations, based on different combinations of spectral complexity coefficients, source locations and $\textup{SNR}^\mbf{X}$ levels. For each configuration we collected the two optimal parameters $\lambda_\mbf{x}^*$ and $\lambda_\mbf{S}^*$ and we investigated their dependence on the signal-to-noise-ratio. In accordance with classical results from inverse theory \citep{hanke1993regularization}, we found that $\lambda_\mbf{x}^*$ depends on the signal-to-noise ratio. What is novel here is the relation between $\lambda_\mbf{S}^*$ and both $\textup{SNR}^\mbf{X}$ and $\textup{SNR}^\mbf{S}$. Indeed for increasing $\textup{SNR}^\mbf{X}$ less regularization is needed, but such dependence varies with the MVAR models. On the other side, the dependence of $\lambda_\mbf{S}^*$ on $\textup{SNR}^\mbf{S}$ is neater and does not depend on the models. Figure \ref{fig:lamS_SNR_S} shows this result; on the left the regularization parameters for the cross-power spectrum reconstruction versus $\textup{SNR}^\mbf{X}$ are shown, while on the right the same parameters are shown with respect to $\textup{SNR}^\mbf{S}$. For the ease of presentation the figure shows the parameters related to one source location; while on the left lines corresponding to different MVAR models have different heights, on the right they overlap. 
\subsection{$ \lambda^*_\mbf{S}<\lambda^*_\mbf{x}$ and dependency from the spectral complexity}

\begin{figure}
    \centering
    \includegraphics[scale=0.3]{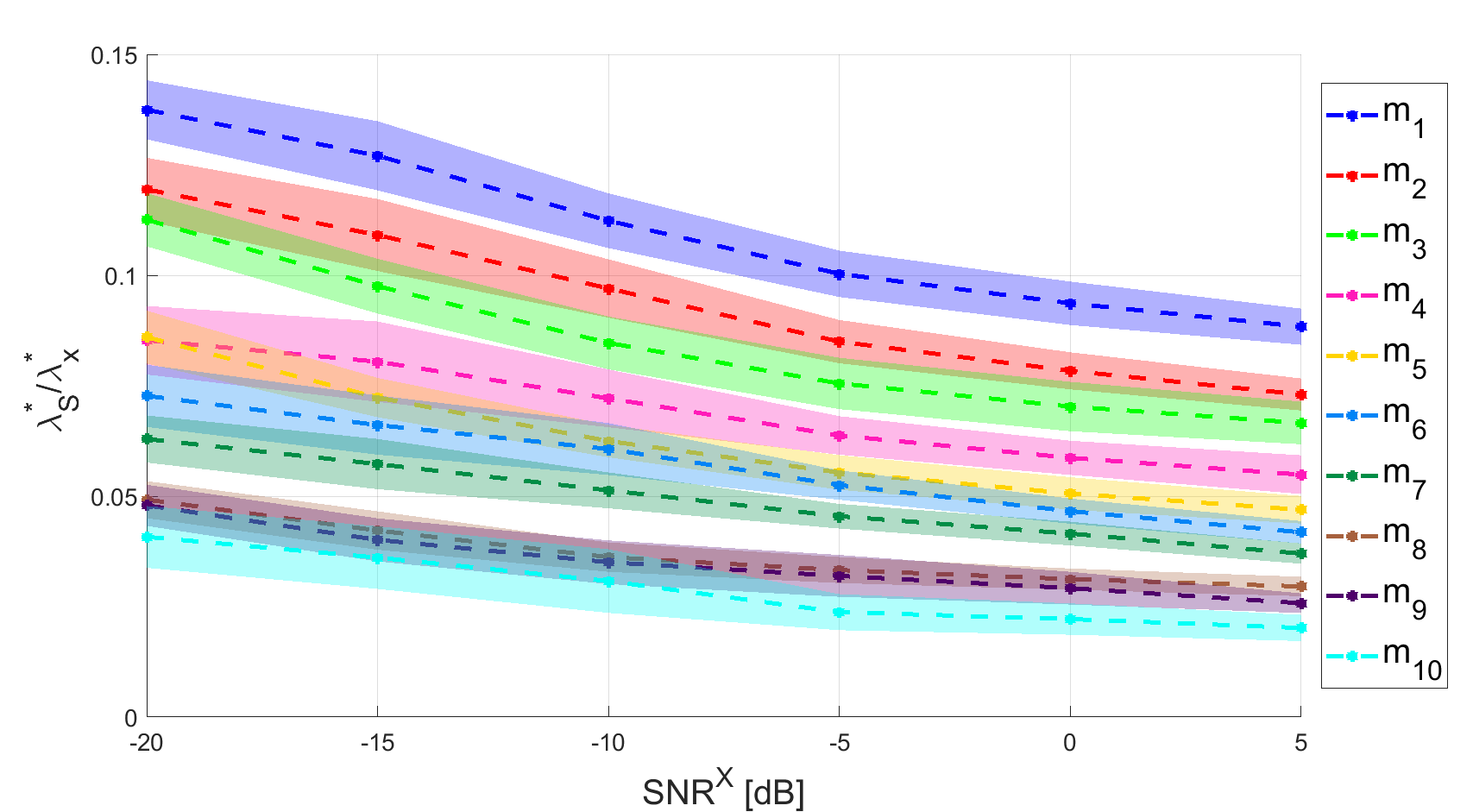}
    \caption{Ratio between the optimal parameters ($\frac{\lambda_\mbf{S}^*}{\lambda_\mbf{x}^*}$) as a function of $\textup{SNR}^\mbf{X}$. Different colours correspond to MVAR model with different spectral complexities. Dashed lines are the mean of the ratio over the different sources location, while solid colours correspond to the standard deviation of the mean.}
    \label{fig:ratio}
\end{figure}

We also investigated the relation between the two optimal regularization parameters. Figure \ref{fig:ratio} shows the ratio $\frac{\lambda_\mbf{S}^*}{\lambda_\mbf{X}^*}$ versus $\textup{SNR}^\mbf{X}$ for the simulated MVAR models. The ratio between the two parameters is always smaller than $\frac{1}{2}$, meaning that $\lambda_\mbf{S}^*<\frac{1}{2}\lambda_\mbf{x}^*$, as it was analytically proved in a simplified case in \cite{vallarino2020two}. Further to this, the figure shows that for increasing spectral complexity coefficients this ratio gets smaller. This latter result is directly related to equation (\ref{eq:SNR_relation}). In fact, for same levels of $\textup{SNR}^\mbf{X}$, signals with higher spectral complexity have higher $\textup{SNR}^\mbf{S}$ and, thus, need less regularization.\\

\section{Discussion and conclusions}
In the present work, we investigated the role of the spectral complexity of a time series, $\mbf{x}(t)$, in the design of an optimal inverse technique for estimating its cross-power spectrum, $\mbf{S}^{\mbf{x}}(f)$, from indirect measurements of the time series itself. Motivated by an analysis pipeline widely used for estimating brain functional connectivity from MEG data, we reconstructed the cross-power spectrum in two steps: first, we estimated the unknown time series by using the Tikhonov method; then we computed the cross-power spectrum of the reconstructed time series. In the present work, we used numerical simulations to study how the spectral complexity of $\mbf{x}(t)$ impacts the value of the regularization parameter that provides the best reconstruction of the cross-power spectrum.

As a first analytical result, we related $\textup{SNR}^\mbf{X}$ to $\textup{SNR}^\mbf{S}$, i.e. the signal-to-noise ratio of the time series and the signal-to-noise ratio of  the corresponding cross-power spectra. The obtained formula suggests that, for a fixed level of $\textup{SNR}^\mbf{X}$,  $\textup{SNR}^\mbf{S}$ depends on the spectral complexity of $\mbf{x}(t)$: the higher the spectral complexity coefficient the higher $\textup{SNR}^\mbf{S}$. Intuitively this happens because a higher value of the spectral complexity coefficient corresponds to a more peaked cross-power spectrum that will emerge over the cross-power spectrum of the noise. 

To test the effect of this result on the choice of the Tikohonov regularization parameter in a practical scenario, we simulated a large set of MEG data and we applied the described two-step approach for estimating the cross-power spectrum of the underlying neural sources. In details, we simulated 1200 synthetic MEG data with varying $\textup{SNR}^\mbf{X}$ generated by pairs of coupled point-like sources at varying locations and with different spectral complexities. For each simulated data we computed the two parameters providing the best estimate of the time series ($\lambda^*_\mbf{x}$) and the best estimate of the cross-power spectrum ($\lambda^*_\mbf{S}$). As shown by Figure \ref{fig:lamS_SNR_S}, the results of our simulations highlighted a high correlation between the values of $\lambda^*_\mbf{S}$ and of $\textup{SNR}^\mbf{S}$.

Eventually, we focused on the relationship between the two parameters $\lambda^*_\mbf{x}$ and $\lambda_\mbf{S}^*$, whose ratio is shown in Figure \ref{fig:ratio}. The figure points out that this ratio depends on the spectral complexity of the simulated time series. This fact may be understood in lights of the previous results, as $\lambda_\mbf{S}^*$ depends on $\textup{SNR}^\mbf{S}$ that in turns depends on the spectral complexity coefficient. Additionally, we found that, for all the simulated data, $\frac{\lambda_\mbf{S}^*}{\lambda^*_\mbf{x}} < \frac{1}{2}$, in line with the results shown in \cite{vallarino2020two} for a simplified model where the neural time series were assumed to be white Gaussian processes. Moreover, when the spectral complexity coefficient increases ($c > 5$ in our simulations) the ratio between the two parameters approaches 0.01. This agrees with the results shown in \cite{hietal16} where, by simulating sinusoidal signals, the authors suggested to use for connectivity estimation a parameter of two orders of magnitude lower.   

The present work focuses on the cross-power spectrum as a connectivity metric. Even though the cross-power spectrum is the starting point for the computation of many connectivity metrics it would be interesting to directly investigate the behaviour of the Thikonov regularization parameters when using such metrics. Future works will be devoted to this. It is also worth noticing that the definition of optimality when defining the regularization parameters is not univocal, since many metrics can be used. A common example is the area under the curve (AUC), which is the metric that was used in \cite{hietal16}. The use of different metrics would firstly strengthen our results and would also allow a more straightforward comparison with the results of \cite{hietal16}. Finally, the dependence of $\lambda_\mbf{S}^*$ on $\textup{SNR}^\mbf{S}$ suggests that an analysis of such dependence could be considered for the definition of a rule for choosing $\lambda_\mbf{S}^*$ in practical scenarios.

\vspace{6pt} 


\funding{EV, AS, SS and MP have been partially supported by Gruppo Nazionale per il Calcolo Scientifico.}


\conflictsofinterest{The authors declare no conflict of interest.} 


\reftitle{References}
\externalbibliography{yes}
\bibliography{main.bib}

\begin{thebibliography}{-------}
\providecommand{\natexlab}[1]{#1}

\bibitem[H{\"a}m{\"a}l{\"a}inen \em{et~al.}(1993)H{\"a}m{\"a}l{\"a}inen, Hari,
  Ilmoniemi, Knuutila, and Lounasmaa]{hamalainen1993magnetoencephalography}
H{\"a}m{\"a}l{\"a}inen, M.; Hari, R.; Ilmoniemi, R.J.; Knuutila, J.; Lounasmaa,
  O.V.
\newblock Magnetoencephalography—theory, instrumentation, and applications to
  noninvasive studies of the working human brain.
\newblock {\em Reviews of modern Physics} {\bf 1993}, {\em 65},~413.

\bibitem[H{\"a}m{\"a}l{\"a}inen and Ilmoniemi(1994)]{hail94}
H{\"a}m{\"a}l{\"a}inen, M.S.; Ilmoniemi, R.J.
\newblock Interpreting magnetic fields of the brain: minimum norm estimates.
\newblock {\em Medical \& biological engineering \& computing} {\bf 1994}, {\em
  32},~35--42.

\bibitem[Calvetti \em{et~al.}(2015)Calvetti, Pascarella, Pitolli, Somersalo,
  and Vantaggi]{calvetti2015hierarchical}
Calvetti, D.; Pascarella, A.; Pitolli, F.; Somersalo, E.; Vantaggi, B.
\newblock A hierarchical Krylov--Bayes iterative inverse solver for {MEG} with
  physiological preconditioning.
\newblock {\em Inverse Problems} {\bf 2015}, {\em 31},~125005.

\bibitem[Costa \em{et~al.}(2017)Costa, Batatia, Oberlin, D'Giano, and
  Tourneret]{costa2017bayesian}
Costa, F.; Batatia, H.; Oberlin, T.; D'Giano, C.; Tourneret, J.Y.
\newblock Bayesian {EEG} source localization using a structured sparsity prior.
\newblock {\em NeuroImage} {\bf 2017}, {\em 144},~142--152.

\bibitem[Sorrentino and Piana(2017)]{sorrentino2017inverse}
Sorrentino, A.; Piana, M.
\newblock Inverse Modeling for {MEG}/{EEG} data. In {\em Mathematical and
  Theoretical Neuroscience}; Springer, Cham,  2017; pp. 239--253.

\bibitem[Bekhti \em{et~al.}(2018)Bekhti, Lucka, Salmon, and
  Gramfort]{bekhti2018hierarchical}
Bekhti, Y.; Lucka, F.; Salmon, J.; Gramfort, A.
\newblock A hierarchical Bayesian perspective on majorization-minimization for
  non-convex sparse regression: application to {M/EEG} source imaging.
\newblock {\em Inverse Problems} {\bf 2018}, {\em 34},~085010.

\bibitem[Luria \em{et~al.}(2019)Luria, Duran, Visani, Sommariva, Rotondi,
  Sebastiano, Panzica, Piana, and Sorrentino]{luria19}
Luria, G.; Duran, D.; Visani, E.; Sommariva, S.; Rotondi, F.; Sebastiano, D.R.;
  Panzica, F.; Piana, M.; Sorrentino, A.
\newblock Bayesian multi-dipole modelling in the frequency domain.
\newblock {\em Journal of neuroscience methods} {\bf 2019}, {\em 312},~27--36.

\bibitem[Ilmoniemi and Sarvas(2019)]{ilsa19}
Ilmoniemi, R.J.; Sarvas, J.
\newblock {\em {B}rain {S}ignals: {P}hysics and {M}athematics of {MEG} and
  {EEG}}; Mit Press,  2019.

\bibitem[Geweke(1982)]{geweke1982measurement}
Geweke, J.
\newblock Measurement of linear dependence and feedback between multiple time
  series.
\newblock {\em Journal of the American statistical association} {\bf 1982},
  {\em 77},~304--313.

\bibitem[Baccal{\'a} and Sameshima(2001)]{baccala2001partial}
Baccal{\'a}, L.A.; Sameshima, K.
\newblock Partial directed coherence: a new concept in neural structure
  determination.
\newblock {\em Biological cybernetics} {\bf 2001}, {\em 84},~463--474.

\bibitem[Nolte \em{et~al.}(2004)Nolte, Bai, Wheaton, Mari, Vorbach, and
  Hallett]{nolte04}
Nolte, G.; Bai, O.; Wheaton, L.; Mari, Z.; Vorbach, S.; Hallett, M.
\newblock Identifying true brain interaction from {EEG} data using the
  imaginary part of coherency.
\newblock {\em Clinical neurophysiology} {\bf 2004}, {\em 115},~2292--2307.

\bibitem[Pereda \em{et~al.}(2005)Pereda, Quiroga, and Bhattacharya]{pereda_05}
Pereda, E.; Quiroga, R.Q.; Bhattacharya, J.
\newblock Nonlinear multivariate analysis of neurophysiological signals.
\newblock {\em Progress in neurobiology} {\bf 2005}, {\em 77},~1--37.

\bibitem[Sakkalis(2011)]{sakkalis_11}
Sakkalis, V.
\newblock Review of advanced techniques for the estimation of brain
  connectivity measured with {EEG}/{MEG}.
\newblock {\em Computers in biology and medicine} {\bf 2011}, {\em
  41},~1110--1117.

\bibitem[Chella \em{et~al.}(2014)Chella, Marzetti, Pizzella, Zappasodi, and
  Nolte]{chella2014third}
Chella, F.; Marzetti, L.; Pizzella, V.; Zappasodi, F.; Nolte, G.
\newblock Third order spectral analysis robust to mixing artifacts for mapping
  cross-frequency interactions in {EEG}/{MEG}.
\newblock {\em Neuroimage} {\bf 2014}, {\em 91},~146--161.

\bibitem[Schoffelen and Gross(2009)]{schoffelen2009source}
Schoffelen, J.M.; Gross, J.
\newblock Source connectivity analysis with MEG and EEG.
\newblock {\em Human brain mapping} {\bf 2009}, {\em 30},~1857--1865.

\bibitem[Barzegaran and Knyazeva(2017)]{barzegaran2017functional}
Barzegaran, E.; Knyazeva, M.G.
\newblock Functional connectivity analysis in EEG source space: the choice of
  method.
\newblock {\em PloS one} {\bf 2017}, {\em 12},~e0181105.

\bibitem[Van~de Steen \em{et~al.}(2019)Van~de Steen, Faes, Karahan, Songsiri,
  Valdes-Sosa, and Marinazzo]{van2019critical}
Van~de Steen, F.; Faes, L.; Karahan, E.; Songsiri, J.; Valdes-Sosa, P.A.;
  Marinazzo, D.
\newblock Critical comments on EEG sensor space dynamical connectivity
  analysis.
\newblock {\em Brain topography} {\bf 2019}, {\em 32},~643--654.

\bibitem[Fries(2005)]{fries2005}
Fries, P.
\newblock A mechanism for cognitive dynamics: neuronal communication through
  neuronal coherence.
\newblock {\em Trends in cognitive sciences} {\bf 2005}, {\em 9},~474--480.

\bibitem[Fries(2015)]{fries2015}
Fries, P.
\newblock Rhythms for cognition: communication through coherence.
\newblock {\em Neuron} {\bf 2015}, {\em 88},~220--235.

\bibitem[Nunez \em{et~al.}(1999)Nunez, Silberstein, Shi, Carpenter, Srinivasan,
  Tucker, Doran, Cadusch, and Wijesinghe]{nunez99}
Nunez, P.L.; Silberstein, R.B.; Shi, Z.; Carpenter, M.R.; Srinivasan, R.;
  Tucker, D.M.; Doran, S.M.; Cadusch, P.J.; Wijesinghe, R.S.
\newblock {EEG} coherency {II}: experimental comparisons of multiple measures.
\newblock {\em Clinical Neurophysiology} {\bf 1999}, {\em 110},~469--486.

\bibitem[Nolte \em{et~al.}(2008)Nolte, Ziehe, Nikulin, Schl{\"o}gl, Kr{\"a}mer,
  Brismar, and M{\"u}ller]{nolte08}
Nolte, G.; Ziehe, A.; Nikulin, V.V.; Schl{\"o}gl, A.; Kr{\"a}mer, N.; Brismar,
  T.; M{\"u}ller, K.R.
\newblock Robustly estimating the flow direction of information in complex
  physical systems.
\newblock {\em Physical review letters} {\bf 2008}, {\em 100},~234101.

\bibitem[Schoffelen and Gross(2019)]{sc_gr19}
Schoffelen, J.M.; Gross, J.
\newblock Studying dynamic neural interactions with {MEG}.
\newblock {\em Magnetoencephalography: from signals to dynamic cortical
  networks} {\bf 2019}, pp. 1--23.

\bibitem[Hincapi{\'e} \em{et~al.}(2016)Hincapi{\'e}, Kujala, Mattout,
  Daligault, Delpuech, Mery, Cosmelli, and Jerbi]{hietal16}
Hincapi{\'e}, A.S.; Kujala, J.; Mattout, J.; Daligault, S.; Delpuech, C.; Mery,
  D.; Cosmelli, D.; Jerbi, K.
\newblock {MEG} Connectivity and Power Detections with Minimum Norm Estimates
  Require Different Regularization Parameters.
\newblock {\em Computational intelligence and neuroscience} {\bf 2016}, {\em
  2016}.

\bibitem[Vallarino \em{et~al.}(2020)Vallarino, Sommariva, Piana, and
  Sorrentino]{vallarino2020two}
Vallarino, E.; Sommariva, S.; Piana, M.; Sorrentino, A.
\newblock On the two-step estimation of the cross-power spectrum for dynamical
  linear inverse problems.
\newblock {\em Inverse Problems} {\bf 2020}, {\em 36},~045010.

\bibitem[Bendat and Piersol(2011)]{be_pi11}
Bendat, J.S.; Piersol, A.G.
\newblock {\em Random data: analysis and measurement procedures}; Vol. 729,
  John Wiley \& Sons,  2011.

\bibitem[Welch(1967)]{welch1967use}
Welch, P.
\newblock The use of fast Fourier transform for the estimation of power
  spectra: a method based on time averaging over short, modified periodograms.
\newblock {\em IEEE Transactions on audio and electroacoustics} {\bf 1967},
  {\em 15},~70--73.

\bibitem[Tikhonov \em{et~al.}(2013)Tikhonov, Goncharsky, Stepanov, and
  Yagola]{ti_etal13}
Tikhonov, A.N.; Goncharsky, A.; Stepanov, V.; Yagola, A.G.
\newblock {\em Numerical methods for the solution of ill-posed problems}; Vol.
  328, Springer Science \& Business Media,  2013.

\bibitem[Chella \em{et~al.}(2019)Chella, Marzetti, Stenroos, Parkkonen,
  Ilmoniemi, Romani, and Pizzella]{chella2019impact}
Chella, F.; Marzetti, L.; Stenroos, M.; Parkkonen, L.; Ilmoniemi, R.J.; Romani,
  G.L.; Pizzella, V.
\newblock The impact of improved MEG--MRI co-registration on MEG connectivity
  analysis.
\newblock {\em Neuroimage} {\bf 2019}, {\em 197},~354--367.

\bibitem[Thompson \em{et~al.}(1991)Thompson, Brown, Kay, and
  Titterington]{thompson1991study}
Thompson, A.M.; Brown, J.C.; Kay, J.W.; Titterington, D.M.
\newblock A study of methods of choosing the smoothing parameter in image
  restoration by regularization.
\newblock {\em IEEE Transactions on Pattern Analysis \& Machine Intelligence}
  {\bf 1991}, pp. 326--339.

\bibitem[Hanke and Hansen(1993)]{hanke1993regularization}
Hanke, M.; Hansen, P.C.
\newblock Regularization methods for large-scale problems.
\newblock {\em Surv. Math. Ind} {\bf 1993}, {\em 3},~253--315.

\bibitem[Hansen(1998)]{hansen1998rank}
Hansen, P.C.
\newblock {\em Rank-deficient and discrete ill-posed problems: numerical
  aspects of linear inversion}; SIAM,  1998.

\bibitem[Vogel(2002)]{vogel2002computational}
Vogel, C.R.
\newblock {\em Computational methods for inverse problems}; SIAM,  2002.

\bibitem[Baillet \em{et~al.}(2001)Baillet, Mosher, and Leahy]{baillet01}
Baillet, S.; Mosher, J.C.; Leahy, R.M.
\newblock Electromagnetic brain mapping.
\newblock {\em IEEE Signal processing magazine} {\bf 2001}, {\em 18},~14--30.

\bibitem[Sorrentino(2010)]{sorrentino10}
Sorrentino, A.
\newblock Particle filters for magnetoencephalography.
\newblock {\em Archives of Computational Methods in Engineering} {\bf 2010},
  {\em 17},~213--251.

\bibitem[Dale and Sereno(1993)]{dale93}
Dale, A.M.; Sereno, M.I.
\newblock Improved localization of cortical activity by combining EEG and MEG
  with MRI cortical surface reconstruction: a linear approach.
\newblock {\em Journal of cognitive neuroscience} {\bf 1993}, {\em
  5},~162--176.

\bibitem[Gramfort \em{et~al.}(2014)Gramfort, Luessi, Larson, Engemann,
  Strohmeier, Brodbeck, Parkkonen, and H{\"a}m{\"a}l{\"a}inen]{gr_etal14}
Gramfort, A.; Luessi, M.; Larson, E.; Engemann, D.A.; Strohmeier, D.; Brodbeck,
  C.; Parkkonen, L.; H{\"a}m{\"a}l{\"a}inen, M.S.
\newblock {MNE} software for processing {MEG} and {EEG} data.
\newblock {\em Neuroimage} {\bf 2014}, {\em 86},~446--460.

\bibitem[L{\"u}tkepohl(2005)]{lutkepohl2005new}
L{\"u}tkepohl, H.
\newblock {\em New introduction to multiple time series analysis}; Springer
  Science \& Business Media,  2005.

\bibitem[Sommariva \em{et~al.}(2019)Sommariva, Sorrentino, Piana, Pizzella, and
  Marzetti]{soetal19}
Sommariva, S.; Sorrentino, A.; Piana, M.; Pizzella, V.; Marzetti, L.
\newblock A comparative study of the robustness of frequency-domain
  connectivity measures to finite data length.
\newblock {\em Brain topography} {\bf 2019}, {\em 32},~675--695.

\end{thebibliography}

\end{document}